\documentclass[10pt,leqno]{article}
\usepackage{amsmath}
\usepackage{amsfonts}
\usepackage{amsthm}
\usepackage{indentfirst}

\usepackage{verbatim}
\usepackage{path}
\usepackage{hyperref}
\usepackage{fancyhdr}
\newtheorem*{Assumption}{Assumption}
\newtheorem{Conjecture}[]{Conjecture}
\newtheorem*{Goal}{Goal}
\newtheorem{Question}[]{Question}
\newtheorem*{Remark}{Remark}

\newcommand{\odip}[2]{o _{#1}\!\left(#2\right)\mathchoice{\!}{}{}{}}
\newcommand{\odi}[1]{\odip{}{#1}}
\newcommand{\Odi}[1]{\mathcal{O}\left(#1\right)}
\newcommand{\Odim}[1]{\mathcal{O}\bigl(#1\bigr)}
\newcommand{\eps}{\varepsilon}
\newcommand{\dx}{\mathrm{d}}
\newcommand{\ii}{\mathrm{i}}
\newcommand{\R}{\mathbb{R}}
\newcommand{\li}{\operatorname{li}}
\newcommand{\eqdef}{\mathchoice{{\,\overset{\mathrm{def}}=}\,}{:=}{:=}{:=}}

\theoremstyle{plain}
\newtheorem{theorem}{\indent\rm T\,h\,e\,o\,r\,e\,m\;}[section]

\theoremstyle{definition}

\theoremstyle{remark}

\makeatletter                                                  %
\renewcommand*{\@seccntformat}[1]{
  \csname the#1\endcsname\;-                                   %
}                                                              %
\renewcommand{\section}{\@startsection{section}{1}{0mm}        %
   {1.5\baselineskip}
   {1\baselineskip}
   {\indent\normalfont\normalsize\bfseries}
   }                                                           %
\renewcommand*{\@seccntformat}[1]{
  \normalfont\bfseries\csname the#1\endcsname\;-               %
}                                                              %
\renewcommand\subsection{\@startsection                        %
  {subsection}{2}{0mm}
  {1.5\baselineskip}
  {1\baselineskip}
  {\indent\normalfont\normalsize\itshape}}
\renewcommand*{\@seccntformat}[1]{
  \normalfont\bfseries\csname the#1\endcsname\;-               %
}                                                              %
\renewcommand\subsubsection{\@startsection                     %
  {subsubsection}{2}{0mm}
  {1.5\baselineskip}
  {1\baselineskip}
  {\indent\normalfont\normalsize\texttt}}
\makeatother                                                   %
\begin{document}
\thispagestyle{empty}

\vskip -8in
\begin{center}
 \rule{8.5cm}{0.5pt}\\[-0.1cm] {\small Riv.\, Mat.\, Univ.\, Parma,\,
Vol. {\bf x} \,(2016), \,000-000}\\[-0.25cm] \rule{8.5cm}{0.5pt}
\end{center}
\vspace {2.2cm}

\begin{center}
{\sc\large Alessandro Zaccagnini}
\end{center}
\vspace {1.5cm}

\centerline{\large{\textbf{The Selberg integral and a new pair-correlation
       function}}}
\smallskip
\centerline{\large{\textbf{for the zeros of the Riemann zeta-function}}}

\renewcommand{\thefootnote}{\fnsymbol{footnote}}

\footnotetext{
This research was partially supported by the grant PRIN2010-11
\textsl{Arithmetic Algebraic Geometry and Number Theory}.}

\renewcommand{\thefootnote}{\arabic{footnote}}
\setcounter{footnote}{0}

\vspace{1,5cm}
\begin{center}
\begin{minipage}[t]{10cm}

\small{ \noindent \textbf{Abstract.}
The present paper is a report on joint work with Alessandro Languasco
and Alberto Perelli, collected in \cite{LanguascoPZ2012a},
\cite{LanguascoPZ2013a} and \cite{LanguascoPZ2013b}, on our recent
investigations on the Selberg integral and its connections to
Montgomery's pair-correlation function.
We introduce a more general form of the Selberg integral and connect
it to a new pair-correlation function, emphasising its relations to
the distribution of prime numbers in short intervals.
\medskip

\medskip

\noindent \textbf{Mathematics~Subject~Classification~(2010):}
11M26, 11N05.
}
\end{minipage}
\end{center}

\bigskip


\section{Introduction}

The central problem of Analytic Number Theory is the distribution of
prime numbers.
How many prime numbers are there, and are they ``randomly'' or
``evenly'' distributed?
The answers to these questions are known only partially, even if one
assumes powerful and, as yet, unproved hypotheses like Riemann's.
Here we are interested in the distribution of prime numbers in ``short
intervals,'' from various points of view.
We give a strong quantitative version of previous results on the
connection between the Selberg integral $J$ defined in \eqref{def-J}
and Montgomery's pair-correlation function $F$ defined in
\eqref{def-F}.
Then we introduce a new, generalised pair-correlation function and
show that it essentially is in charge of the distribution of primes in
short intervals, giving a unified view of many problems.

While it is always advisable to read the original papers (like
Montgomery's \cite{Montgomery1973} and the other ones quoted at
appropriate places), it is extremely instructive to start with
Goldston's excellent survey \cite{Goldston2005}.

\section{Distribution of prime numbers: Large scale}

\subsection{The Chebyshev prime-counting functions}

Let
\[
  \pi(x)
  \eqdef
  \vert \{ p \le x \colon p \text{ is prime} \} \vert
  \qquad\text{and}\qquad
  \psi(x)
  \eqdef
  \sum_{n \le x} \Lambda(n),
\]
where $\Lambda(n) = 0$ unless $n = p^\alpha$ for some prime $p$ and
positive integer $\alpha$; and $\Lambda \bigl( p^\alpha \bigr) = \log p$.
We also let
\[
  \li(x)
  \eqdef
  \int_2^x \frac{\dx t}{\log t}.
\]
It has been conjectured by Gauss at the end of the XVIII century that
$\pi(x) \sim \li(x)$ as $x \to +\infty$.
Stated in this way, Gauss's conjecture is equivalent to Legendre's
simpler statement that $\pi(x) \sim x / \log x$, but Gauss's
approximation is numerically far closer to the truth.
The truth of the asymptotic statement above was proved independently
by Hadamard and de la Vall\'ee Poussin in 1896 and is called the Prime
Number Theorem (PNT, for short).
Let
\[
  \Delta_\pi(x)
  \eqdef
  \pi(x)
  -
  \li(x)
  \qquad\text{and}\qquad
  \Delta_\psi(x)
  \eqdef
  \psi(x)
  -
  x.
\]
According to the strongest form of the PNT known today, there is a
constant $c > 0$ such that for $x \to +\infty$
\begin{equation}
\label{PNT}
  \Delta_\pi(x),
  \Delta_\psi(x)
  \ll
  R(x)
  \eqdef
  x \exp \left\{ -c (\log x)^{3/5} (\log\log x)^{-1/5} \right\}.
\end{equation}
The PNT shows that to a first approximation $\pi(x)$ is very close to
$\li(x)$ and $\psi(x)$ is very close to $x$.
Of course, it is interesting, and in some applications critical, to
know the exact degree of these approximations.
It has been known since Riemann's pioneering paper \cite{Riemann1859}
that $\Delta_\pi$ and $\Delta_\psi$ both depend on the distribution of
the zeros of the zeta-function.
The connection has received a strong quantitative form by Pintz
\cite{Pintz1984}, building on earlier work started by Ingham.
Here we are mainly concerned with the case of ``optimal
distribution,'' and we will not pursue this connection further.

\subsection{The Riemann Hypothesis}

In his 1859 paper referred to above, Riemann proved that the
zeta-function defined by the Dirichlet series
\[
  \zeta(s)
  =
  \sum_{n \ge 1} \frac1{n^s}
  \qquad\text{for $\sigma = \Re(s) > 1$,}
\]
has a meromorphic continuation to the complex plane, except for the
point $s = 1$ where $\zeta$ has a simple pole with residue $1$.
The continuation has infinitely many real zeros (the so-called
``trivial zeros'') at the points $-2$, $-4$, $-6$, \dots, and also
infinitely many non-real zeros in the strip $\sigma \in [0, 1]$.
It is customary to denote the generic non-trivial zero of the Riemann
$\zeta$-function by $\rho \eqdef \beta + \ii \gamma$.
These zeros are placed symmetrically with respect to the real axis and
the line $\sigma = \frac12$.
Riemann wrote that it looks ``likely'' that all of these zeros actually
lie on the line $\sigma = \frac12$, but that he was unable to prove
it.
This is the original statement of the Riemann Hypothesis (RH, for short)
that is still unsettled, as yet.

The Riemann Hypothesis is equivalent to either of the two statements
\begin{equation}
\label{RH}
  \pi(x)
  =
  \int_2^x \frac{\dx t}{\log t}
  +
  \Odi{ x^{1/2} \log x }
  \quad\text{or}\quad
  \psi(x)
  =
  x
  +
  \Odi{ x^{1/2} (\log x)^2 }.
\end{equation}
Of course, this is \emph{much stronger} than \eqref{PNT}.
Notice that by Littlewood's Theorem \cite{Littlewood1914} of 1914, we
have
\begin{equation}
\label{Delta-JEL}
  \Delta_\psi(x)
  =
  \psi(x) - x
  =
  \Omega_{\pm}\bigl( x^{1 /2} \log \log \log x \bigr).
\end{equation}
Here, the notation means that $\Delta_\psi(x)$ is larger than a positive
multiple of $x^{1 /2} \log \log \log x$ for a suitable unbounded
sequence of values of $x$, and smaller than a negative multiple of the
same function on another such sequence.
In other words, the RH is very nearly optimal.
It seems plausible that the correct maximal order of magnitude
for $\psi(x) - x$ be $x^{1 / 2} (\log \log \log x)^2$.
See \S15.3 of Montgomery \& Vaughan \cite{MontgomeryV2007}.

\section{Distribution of prime numbers: Small scale}

\subsection{Primes in all ``short'' intervals}

The ``additive'' form of the expected main term for both $\pi$ and
$\psi$ readily suggests a conjecture on the number of prime numbers
contained in the interval $(x, x + y]$.
We say that such an interval is \emph{short} if $y = o(x)$.

\begin{Question}
For $y \le x$, is it true that
\begin{equation}
\label{asymp-Delta}
  \pi(x + y) - \pi(x)
  \sim
  \int_x^{x + y} \frac{\dx t}{\log t}
  \quad\text{or}\quad
  \psi(x + y) - \psi(x)
  \sim
  y \quad ?
\end{equation}
\end{Question}

The asymptotic relations in \eqref{asymp-Delta} are true if
$y / R(x) \to + \infty$.
This is essentially trivial, because we can use \eqref{PNT} twice and
the difference between the leading terms is larger than the error
terms.

We may rephrase the question as: Is \eqref{asymp-Delta} true also for
smaller $y$?
For \emph{very} small $y$ (say, $y \asymp \log x$ or smaller) it is
trivially false, because most intervals will not be long enough to
contain any primes at all.
On the other hand, if $y = 1$ and $x + 1$ is a prime number, then
the interval $(x, x + 1]$ contains ``too many'' primes.
We assume that $y \to + \infty$ with $x$ in order to avoid such
trivialities.

The asymptotic formulae \eqref{asymp-Delta} are false also for
$y = (\log x)^\alpha$, for any \emph{fixed} $\alpha > 1$: this has
been proved by Maier in 1985 \cite{Maier1985}, and it came as a
surprise since, as we know conditionally from Selberg's work
\cite{Selberg1943}, they are \emph{usually} true; see the discussion
below.
Unconditionally, the relations \eqref{asymp-Delta} are true for
$y \ge x^{7 / 12 -\eps(x)}$, provided that $\eps(x) \to 0+$ as
$x \to +\infty$; see Heath-Brown \cite{HeathBrown1988}.
Of course, on the RH they are true in a wider range: by \eqref{RH},
for $y / \bigl( x^{1 / 2} (\log x)^2 \bigr) \to + \infty$, but we must
recall that Cram\'er proved that the weaker condition
$y / \bigl( x^{1 / 2} (\log x) \bigr) \to + \infty$ suffices.
Actually, his techniques also prove that the RH implies the existence
of at least $x^{1 / 2}$ prime numbers in the interval $[x, x + y]$
provided that $y \ge C x^{1 / 2} \log x$, where $C$ is a sufficiently
large constant.
See Theorem~13.3 in \cite{MontgomeryV2007}.

A long-standing open problem is to bridge the gap between Maier-type
results and what is known under the RH.
It is plausible that we can take $y = x^\eps$ for any fixed $\eps > 0$,
or even slightly smaller.
As we remarked above, by Littlewood's Theorem \eqref{Delta-JEL} we
know that the error term $\Delta_\psi(x) = \psi(x) - x$ is sometimes
quite large, but it is conceivable that it varies very slowly and then
one can compute accurately $\psi(x + y) - \psi(x)$ for comparatively
small $y$, although indirectly.

These topics are discussed at great length and in full detail in
Chapters 13 and 15 of Montgomery \& Vaughan \cite{MontgomeryV2007}.
See, in particular, \S13.1.

\subsection{Primes in ``almost all'' short intervals: The Selberg integral}
\label{sec:primes-aa}

In some applications, it is not necessary to know that either of
\eqref{asymp-Delta} holds \emph{for every} $y$, but only that it
holds for most values of $y$.
Selberg introduced a very convenient way of measuring precisely what
``usually'' means.
Let
\begin{equation}
\label{def-J}
  J(x, \theta)
  \eqdef
  \int_x^{2 x}
    \vert \psi(t + \theta t) - \psi(t) - \theta t \vert^2 \, \dx t
\end{equation}
denote the ``variance'' of the primes in short intervals, where, in
the notation of the previous sections, $\theta \in [0, 1]$ is
essentially $y / x$.

\begin{Question}
In which range of values for $\theta$ is it true that
\begin{equation}
\label{bound-J}
  J(x, \theta)
  =
  \odi{x^3 \theta^2} \quad ?
\end{equation}
\end{Question}

If $\theta$ is not too small the Brun-Titchmarsh inequality in Theorem
\ref{BT} implies that $J(x, \theta) \ll x^3 \theta^2$.
Furthermore, we know that \eqref{bound-J} holds unconditionally
in the range $x^{-5/6 - \eps(x)} \le \theta \le 1$, provided that
$\eps(x) \to 0$ as $x \to \infty$; see \cite{Zaccagnini1998}.
A lower bound for $\theta$ is needed, since for $\theta = 1 / x$, say,
$J(x, \theta)$ essentially reduces to
\begin{equation}
\label{theta-small}
  \sum_{n \in [x, 2 x]} (\Lambda(n) - 1)^2
  \sim
  x \log x,
\end{equation}
by the Prime Number Theorem, and hence \eqref{bound-J} is false in
this range.
Actually, the Brun-Titchmarsh inequality in Theorem \ref{BT} implies,
more generally, that
$J(x, \theta) \ll x^3 \theta^2 (\log x)^2 (\log (2 \theta x))^{-2}$ for
$\theta \ge x^{-1}$, which is compatible with \eqref{theta-small}.
Notice, however, that the mentioned inequality sometimes ``loses'' a
log-factor.

If we assume the RH we have a much stronger result, which is due to
Selberg in \cite{Selberg1943}:
\begin{equation}
\label{J-RH}
  J(x, \theta)
  \ll
  x^2 \theta \bigl( \log (2 / \theta) \bigr)^2
  \qquad\text{uniformly for $x^{- 1} \le \theta \le 1$.}
\end{equation}
%
In other words, in this range of values for $\theta$,
\[
  \psi(t + \theta t) - \psi(t)
  =
  \theta t + \Odi{r(x, \theta)}
  \qquad\text{for ``almost all'' $t \in [x, 2 x]$,}
\]
provided that
$r(x, \theta) / \bigl( (\theta x)^{1/2} \log x \bigr) \to +\infty$.
More precisely, the cardinality of the set of integers $t \in [x, 2 x]$
such that the asymptotic formula above fails to hold is $\odi{x}$ as
$x \to +\infty$.
When $\theta = x^{-1} (\log x)^\alpha$ we still have that the expected
asymptotic formula holds for most integers: Maier's result implies
that, in fact, there are infinitely many exceptions, although they
are rather sparse.

Notice that for $\theta = 1 / x$ this strong result is essentially
empty!
Notice also what happens for $\theta = 1$: the average bound is
slightly stronger than the pointwise one given in \eqref{RH}.

The proof of Selberg's results begins observing that, without any
hypothesis, we have
\[
  J(x, \theta)
  \ll
  \sum_{\rho_1}
    \sum_{\rho_2}
      x^{1 + 2 \beta_1}
      \frac{\min(\theta^2, \gamma_1^{-2})}{(1 + | \gamma_1 - \gamma_2|)^2}.
\]
If the RH holds, then $\rho = \frac12 + \ii \gamma$ for all zeros
of the Riemann zeta-function, and therefore
\[
  J(x, \theta)
  \ll
  x^2
  \sum_{\gamma > 0}
    \min(\theta^2, \gamma^{-2}) \log(\gamma),
\]
by the Riemann-von Mangoldt formula \ref{R-vM}.
Considering separately the ranges $\gamma \in (0, \theta^{-1})$ and
$\gamma > \theta^{-1}$, we see that the Riemann Hypothesis implies
\eqref{J-RH}.
Assuming a less restrictive hypothesis on the distribution of the
zeros, one obtains the result described in the next section.

\subsection{Strong bounds for \texorpdfstring{$J$}{J} and the quasi-RH}

Recall that $\rho \eqdef \beta + \ii \gamma$ denotes the generic
non-trivial zero of the Riemann $\zeta$-function and let
\[
  \Theta
  \eqdef
  \sup \left\{ \beta \le 1 \colon \zeta(\beta + \ii \gamma) = 0
       \right\}.
\]
Essentially, $\Theta < 1$ implies that
%
$
  J(x, \theta)
  \ll
  x (\theta x)^{2 \Theta + \eps}$
%
for any $\eps > 0$, uniformly for $x^{\eps - 1} \le \theta \le 1$.
On the RH we have $\Theta = \frac12$ and we have the stronger
bound given in \eqref{J-RH};
see e.~g., Saffari \& Vaughan \cite{SaffariV1977}, Lemma 6.
We remark that we do not even know whether $\Theta < 1$ or not, as yet.

It is interesting that it is also possible to invert the arrow, as in
Pintz's results quoted above.
We can deduce information concerning the horizontal distribution of
the zeros of the Riemann zeta-function from hypothetical strong upper
bounds for $J$.

\begin{Remark}
Assume that $J(x, 1) \ll x^\delta$ for some $\delta \in [2, 3]$.
Then $\Theta \le \frac12 (\delta - 1)$
\end{Remark}

Notice that \emph{no} uniformity is needed.
This is a result in the author's paper \cite{Zaccagnini2000b}.
The proof is little more than Cauchy's inequality.

\subsection{Application to an approximation problem}

The results described in the previous section are obviously
interesting in themselves.
Here we quickly quote an application to an approximation problem which
is relevant for the circle method.
We have two exponential sums, one containing explicitly the prime
numbers, and we know that they are close to each other in the
neighbourhood of $0$.
The Selberg integral gives, essentially, the $L^2$-norm of the
difference.
Let
\[
  S(\alpha)
  \eqdef
  \sum_{n \le x} \Lambda(n) \mathrm{e}^{2 \pi \ii n \alpha}
  \qquad\text{and}\qquad
  T(\alpha)
  \eqdef
  \sum_{n \le x} \mathrm{e}^{2 \pi \ii n \alpha}.
\]
It is well known that $T$ gives a good approximation to $S$ for
$\alpha$ small since, in fact, $S(0) = \psi(x)$ and
$T(0) = \lfloor x \rfloor$.
In Diophantine problems like the one described by Br\"udern, Cook \&
Perelli in \cite{BrudernCP1997}, the circle is replaced by the whole
real line which is dissected into three pieces: the main term arises
as an integral over a neighbourhood of $0$ and the width of the
neighbourhood is reflected in the quality of the final result.
According to Lemma~1 in \cite{BrudernCP1997} we have
\[
  \int_{- 1 / (\theta x)}^{1 / (\theta x)}
    \vert S(\alpha) - T(\alpha) \vert^2 \, \dx \alpha
  \ll
  \frac1{(\theta x)^2}
  \bigl( J(x, \theta) + x \bigr)
  +
  \theta x (\log x)^2.
\]
Here we see the important role of uniformity in $\theta$:
the length of the integration interval depends on it.
Ideally, we would like to have $\theta = 1 / (2 x)$ so that there is
no ``minor arc,'' which, by definition, is the set of $\alpha$'s where
we have comparatively poor information on the exponential sum $S$.

\section{Montgomery's pair-correlation function}

We now assume RH until the end.
Following Montgomery, we let
\begin{equation}
\label{def-F}
  F(x, T)
  \eqdef
  \sum_{\gamma_1, \gamma_2 \in [0, T]}
    \frac{4 x^{\ii (\gamma_1 - \gamma_2)}}{4 + (\gamma_1 - \gamma_2)^2}
  =
  \sum_{\gamma_1, \gamma_2 \in [0, T]}
    x^{\ii (\gamma_1 - \gamma_2)} \omega(\gamma_1 - \gamma_2),
\end{equation}
where $\omega(x) = 4 / (4 + x^2)$.
The weight $\omega$ arises naturally as the Fourier transform of
$e^{-2 |x|}$, apart from a normalisation factor.
Montgomery \cite{Montgomery1973} proved that
\[
  F(x, T)
  \sim
  \frac T{2 \pi} \log x
  \qquad\text{as $T \to +\infty$}
\]
uniformly for $T^\eps \le x \le T$, and gave the following conjecture.

\begin{Conjecture}[Montgomery]
\label{Mont-PC}
For $T \to +\infty$ we have
\[
  F(x, T)
  \sim
  \frac T{2 \pi} \log T
  \qquad\text{uniformly for $T \le x \le T^A$.}
\]
\end{Conjecture}

In other words, only the ``diagonal'' terms where $\gamma_1 = \gamma_2$
give a contribution, in view of the Riemann-von Mangoldt formula in
Theorem~\ref{R-vM} below, whereas the other terms, more or less,
cancel out.

\subsection{The link between \texorpdfstring{$F$ and $J$}{F and J}}

A number of authors studied the connection between hypothetical
asymptotic formulae for $J$ and $F$, with increasing precision.
The following table lists the major contributions.

\begin{itemize}

\item
Goldston \& Montgomery (1987) \cite{GoldstonM1987}

\item
Montgomery \& Soundararajan (2002) \cite{MontgomeryS2002}

\item
Chan (2003) \cite{Chan2003}

\item
Languasco, Perelli \& Z. (2012) \cite{LanguascoPZ2012a}.

\end{itemize}

Apart from a possible lack of precision, the argument leading to
\eqref{J-RH} sketched at the end of \S\ref{sec:primes-aa} shows that
$J$ depends crucially on the distribution of zeros of $\zeta$:
Injecting the information contained in Montgomery's Conjecture
\ref{Mont-PC}, one gets the one-term asymptotic formula
\begin{equation}
\label{GM-F->J}
  J(x, \theta)
  \sim
  \frac32 x^2 \theta
  \log(1 / \theta),
\end{equation}
and conversely, superseding \eqref{J-RH}.
More precisely, if the asymptotic relation in Conjecture \ref{Mont-PC}
holds uniformly for
$T \in \bigl[ x^{B_1} L^{-3}, x^{B_2} L^3 \bigr]$ then \eqref{GM-F->J}
holds uniformly for $\theta \in \bigl[ x^{-B_2}, x^{-B_1} \bigr]$,
where $L = \log x$ and $0 < B_1 \le B_2 \le 1$ are fixed.
A similar assumption on the range for \eqref{GM-F->J} allows to
recover Conjecture~\ref{Mont-PC} in the corresponding range.
This is Theorem~2 of Goldston \& Montgomery \cite{GoldstonM1987}.
From now on, we drop all reference to uniformity ranges, for
simplicity.
More generally, the authors mentioned above set out to pursue
the following

\begin{Goal}
Compare the size of the ``error terms'' $R_J$ and $R_F$ in
expansions like
\begin{align*}
  J(x, \theta)
  &=
  c_1 x^2 \theta \log(1 / \theta) + c_2 x^2 \theta + R_J(x, \theta) \\
  F(x, T)
  &=
  c_3 T \log T + c_4 T + R_F(x, T)
\end{align*}
in suitable ranges of uniformity.
\end{Goal}

In both cases the main parameter is $x$.
The expected values of the constants $c_j$ are known.
We may therefore restate our goal in a more precise form.
Let $R_J$ and $R_F$ be defined implicitly by
\[
  J(x, \theta)
  =
  \frac32 x^2 \theta
  \bigl(\log(1 / \theta) + 1 - \gamma - \log(2 \pi) \bigr)
  +
  R_J(x, \theta)
\]
and
\[
  F(x, T)
  =
  \frac T{2 \pi} \Bigl( \log \frac T{2 \pi} - 1 \Bigr)
  +
  R_F(x, T).
\]
We expect that $R_J$ is small if and only if $R_F$ is, provided that
the assumed bound holds in a wide enough range of values.
In fact, assuming RH Montgomery \& Soundararajan proved that,
essentially, $R_J(x, \theta) = \odi{x^2 \theta}$
if and only if the imaginary parts of the zeros of the Riemann
zeta-function are ``well distributed'' in the sense that
$R_F(x, T) = \odi{T}$.
In other words, the diagonal contribution $\gamma_1 = \gamma_2$
dominates, as expected.
Actually, it seems necessary to assume slightly more than just
$R_F(x, T) = \odi{T}$ in order to get
$R_J(x, \theta) = \odi{x^2 \theta}$: see the detailed comment to
equation (1.8) in \cite{LanguascoPZ2012a}.
We stress again the importance of the uniformity ranges in the
hypothesis.

Following Chan's results in \cite{Chan2003}, Languasco, Perelli and
the author studied relations between hypothetical bounds of the type
\[
  R_J(x, \theta)
  =
  \Odi{x^2 \theta^{1 + \alpha}}
  \qquad\text{and}\qquad
  R_F(x, T)
  =
  \Odim{T^{1 - \beta}}.
\]
We state a weakened and simplified form of our results: for full
details see \cite{LanguascoPZ2012a}.
Essentially, for $\alpha$, $\beta > 0$, we have
\begin{align}
\label{hypoth-RJ}
  R_J(x, \theta)
  &\ll
  x^2 \theta^{1 + \alpha}
  &&\Longrightarrow&
  R_F(x, T)
  &\ll
  T^{1 - \alpha / (\alpha + 3)} \\
\label{hypoth-RF}
  R_F(x, T)
  &\ll
  T^{1 - \beta}
  &&\Longrightarrow&
  R_J(x, \theta)
  &\ll
  x^2 \theta^{1 + \beta / 2}.
\end{align}
Here we are totally neglecting log-powers and uniformity in the
various parameters, for simplicity of statement, since our results are
as general as they are cumbersome.
According to the heuristics given by Montgomery \& Soundararajan in
\cite{MontgomeryS2002}, a plausible range for $\alpha$ and $\beta$ is
the interval $(0, \frac12)$.

For the proof we need various forms of abelian-tauberian results,
extending the technique introduced by Goldston \& Montgomery in
\cite{GoldstonM1987}.
Essentially, we have a hypothetical average of a certain function, say
$J$, and we need to transform it into a different average of the same
function, before we can use standard Fourier-transform techniques.
This leads to some loss of uniformity.
This also shows why we use $J$ instead of the perhaps more ``natural''
\[
  \widetilde{J}(x, y)
  \eqdef
  \int_x^{2 x}
    \vert \psi(t + y) - \psi(t) - y \vert^2 \, \dx t,
\]
because $J$ is more directly connected to sums over zeros, via the
explicit formula for $\psi$.
In fact, using $J$ we see that each summand in the finite sum over
zeros factors as $t^\rho c(\rho, \theta)$ as in \eqref{Delta^2_psi}
below, and this leads to simpler handling.

However, there is a standard way to relate $J(x, \theta)$ and
$\widetilde{J}(x, \theta x)$ and the two asymptotic formulae differ
only in the ``secondary main term,'' that is, in the value of the
constant $c_2$ above.
For full details, see Saffari \& Vaughan \cite{SaffariV1977}.

We select a few, significant, steps of the proof.
The first one is the approximate equality
\[
  F(x, T)
  =
  \frac1\pi
  \int_0^T
    \Bigl\vert \sum_\gamma \frac{x^{\ii \gamma}}{1 + (t - \gamma)^2}
    \Bigr\vert^2 \, \dx t
  +
  \Odi{(\log T)^3}.
\]
Let $c(\rho, \theta) \eqdef \bigl( (1 + \theta)^\rho - 1 \bigr) / \rho$.
The first application of a tauberian result yields the asymptotic
formula for
\[
  \int_{\R}
    \Bigl\vert
      \sum_{|\gamma| \le Z} c(\rho, \theta) \frac{x^{\ii \gamma}}{1 + (t - \gamma)^2}
    \Bigr\vert^2 \, \dx t,
\]
where the interval $[0, T]$ has been replaced by the whole real line,
only the zeros with ``small'' imaginary parts are included and there
is the coefficient of the ``right'' shape, namely $c$.
Another application coupled with the Plancherel formula gives the
asymptotic formula for
\[
  \int_x^{2 x}
    \Bigl\vert
      \sum_{|\gamma| \le Z} c(\rho, \theta) t^{\rho}
    \Bigr\vert^2 \, \dx t,
\]
where $Z$ is at our disposal.
This is essentially $J$, since
\begin{equation}
\label{Delta^2_psi}
  \psi(t + \theta t) - \psi(t) - \theta t
  =
  -
  \sum_{|\gamma| \le T} c(\rho, \theta) t^\rho
  +
  \Odi{ \frac tT \bigl( \log tT \bigr)^2 },
\end{equation}
by the Explicit formula \ref{EF}, where $T$ is a free parameter.
Recalling that $\beta < 1$ for every zero of the Riemann
$\zeta$-function, we have
\[
  c(\rho, \theta)
  =
  \int_1^{1 + \theta} t^{\rho - 1} \, \dx t
  \ll
  \min\Bigl( \theta, \frac1{|\gamma|} \Bigr).
\]
This means that we have to choose, essentially, $T = \theta^{-1}$ and
a suitable value for $Z$.
Among other things, here we improve a crucial technical lemma (Lemma~1
of \cite{LanguascoPZ2012a}) that allows us to reach the square root in
\eqref{hypoth-RF} instead of the fourth root as in Chan
\cite{Chan2003}.

\subsection{A different version of Montgomery's pair-correlation conjecture}

A slightly weaker form of Montgomery's Conjecture \ref{Mont-PC}
describes the distribution of the ``gaps'' between zeros of the
Riemann zeta-function.
For the detailed relationship between the two Conjectures, see
\cite{Goldston2005}, Theorem~4.

\begin{Conjecture}[Montgomery]
For fixed $\alpha$ and $\beta$ with $0 < \alpha < \beta$ we have
\[
  \sum_{\substack{\gamma_1, \gamma_2 \in [0, T] \\ ((\gamma_1 - \gamma_2) \log T) / 2 \pi \in [\alpha, \beta]}}
    1
  \sim
  \frac1{2 \pi} T \log T
  \int_\alpha^\beta
    \Bigl( 1 - \Bigl( \frac{\sin(\pi u)}{\pi u} \Bigr)^2 \Bigr) \, \dx u,
\]
as $T \to + \infty$.
\end{Conjecture}

The function inside the integral on the far right is known as the
pair-correlation function for the zeros of the Riemann
$\zeta$-function.
In view of the Riemann-von Mangoldt formula \ref{R-vM}, the average
spacing between consecutive zeros of the Riemann zeta-function with
imaginary parts in $[0, T]$ is $2 \pi / \log T$.

\section{A general pair-correlation function}

\subsection{A unifying approach}

Let $\tau \in [0, 1]$ and define
\[
  F(x, T, \tau)
  \eqdef
  \sum_{\gamma_1, \gamma_2 \in [-T, T]}
    \frac{4 x^{\ii (\gamma_1 - \gamma_2)}}{4 + \tau^2 (\gamma_1 - \gamma_2)^2}.
\]
It is clear that $F(x, T, 1)$ is (essentially) the same as $F(x, T)$.
This function is already present, albeit in a slightly different
guise, in Heath-Brown \& Goldston \cite{HeathBrownG1984}; they do not
directly pursue its relations with primes in short intervals, but use
it simply as a technical device.

Moreover, $F(x^{1 / \tau}, T, \tau)$ is the pair-correlation function
for $Z_\tau(s) = \zeta(s / \tau)$, where $Z_\tau$ is (almost) an
element of the Selberg Class of degree and conductor
\[
  d
  =
  \frac1\tau
  \qquad\text{and}\qquad
  q
  =
  \Bigl( \frac1\tau \Bigr)^{1 / \tau}
\]
respectively.
This means that we can use well-known heuristics for the Selberg class
to make plausible guesses concerning $F(x, T, \tau)$.
Of course, $Z_\tau$ is an element of the Selberg Class only for
$\tau = 1$.
See Murty \& Perelli \cite{MurtyP1999} for the discussion of the
pair-correlation function in the Selberg Class.

\subsection{Properties of the general pair-correlation function}

Notice that $F(x, T, 0) = \vert \Sigma(x, T) \vert^2$ where
\begin{equation}
\label{def-Sigma}
  \Sigma(x, T)
  \eqdef
  \sum_{|\gamma| \le T} x^{\ii \gamma}
\end{equation}
is the exponential sum that appears in Landau's explicit formula
\ref{Landau}.
There is a strong conjecture of Gonek \cite{Gonek1993} concerning the
behaviour of $\Sigma$: see Conjecture~\ref{Gonek-C}.
We remark that $F(x, T, \tau)$ is difficult to estimate for $\tau$
very small ($\le 1 / T$); in fact the trivial bound
\[
  F(x, T, \tau)
  \ll
  \min \bigl( T, \tau^{-1} \bigr) T \log^2 T
\]
becomes very large.
When $\tau = 1$, the trivial bound for $F(x, T)$ is only slightly
larger than the expected truth, by just a log factor.

\subsection{The generalised Selberg Integral}

From now on we assume that $\tau > 0$ in order to avoid trivial
statements.
Let
\[
  J(x, \tau, \theta)
  \eqdef
  \int_x^{x (1 + \tau)}
    \vert \psi(t + \theta t) - \psi(t) - \theta t \vert^2 \, \dx t.
\]
Here we are dealing with ``short intervals'' in two different ways.
The obvious conjecture is $J(x, \tau, \theta) \ll x^{2 + \eps} \tau \theta$.

\begin{Assumption}[Hypothesis $H(\eta)$]
We assume that the RH holds and that
\[
  F(x, T, \tau)
  \ll
  T x^{\eps}
  \qquad\text{uniformly for}\qquad
  \begin{cases}
    x^{\eta} \le T \le x \\
    x^{\eta} / T \le \tau \le 1
  \end{cases}
\]
for some fixed $\eta > 0$ and every $\eps > 0$.
\end{Assumption}

This assumption is ``justified'' by Gonek's Conjecture \ref{Gonek-C}
for $\tau$ small, and by an obvious generalisation of Montgomery's for
$\tau$ large.
In \cite{LanguascoPZ2013a}, Languasco, Perelli and the author give a
variant of $H(\eta)$ valid for $\eta = 0$, which we omit for brevity.
The following results are Theorem~1 and 2, respectively, of
\cite{LanguascoPZ2013a}.

\begin{theorem}
\label{Thm1-II}
If assumption $H(\eta)$ holds for some $\eta \in (0, 1)$, then
\[
  J(x, \tau, \theta)
  \ll
  x^{2 + \eps} \tau \theta
\]
uniformly for $x^{-1} \le \theta \le x^{- \eta}$ and
$\theta x^\eta \le \tau \le 1$.
\end{theorem}

As an immediate corollary, we have that
\[
  \psi(x + y) - \psi(x)
  =
  y + \Odim{y^{1 / 2} x^\eps}
\]
for ``almost all'' $x \in [x, x (1 + \tau)]$ and
$y \in [1, x^{1 - \eta}]$.

\begin{theorem}
\label{Thm2-II}
Let $\eps > 0$ be small.
If assumption $H(\eta)$ holds for some $\eta \in (0, 1 / 2 - 5 \eps)$, then
\[
  \psi(x + y)
  -
  \psi(x)
  =
  y
  +
  \begin{cases}
    \Odi{y^{2 / 3} x^{\eta / 3 + \eps}}
    &\text{for $x^{\eta + 5 \eps} \le y \le x^{1 / 2}$,} \\
    \Odi{y^{1 / 3} x^{1 / 6 + \eta / 3 + \eps}}
    &\text{for $x^{1 / 2} \le y \le x^{1 - \eta}$.}
  \end{cases}
\]
\end{theorem}

We quote (2.3) of \cite{LanguascoPZ2013a} to show the explicit
connection between the ``short'' Selberg integral and the generalised
pair-correlation function that we use in the proof of
Theorem~\ref{Thm1-II}.
In fact we have
\[
  \int_x^{x + y}
    \Big \vert \sum_{\vert \gamma \vert \le U}
      \frac{(t + h)^\rho - t^\rho}{\rho} \Big \vert^2 \, \dx t
   \ll
   \frac{h^2 y}{x}
   \max_{x \le t \le x + h} F \Bigl(t, U, \frac yt \Bigr).
\]
The proof of Theorem~\ref{Thm2-II} uses the ``inertia'' property of
$\psi$: $\psi(x + y) - \psi(x)$ does not change by more than
$\asymp \log x$ as $x$ increases by $1$.
Hence, if
\[
  \psi(x + y)
  -
  \psi(x)
  -
  y
\]
is abnormally large in absolute value for some $x = x_0$, it is also
large for ``many'' values of $x$ around $x_0$, which is impossible
by the previous Theorem~\ref{Thm1-II}.
We also adapt an argument in Gonek \cite{Gonek1993}.
In the proof of Theorem~\ref{Thm1-II} we need suitable, uniform
generalisations of the results in the original paper by Heath-Brown
\cite{HeathBrown1982a}.

\subsection{More consequences}

One of the more interesting corollaries of our work in
\cite{LanguascoPZ2013a} is the fact that we obtain explicit estimates
for $\Delta_\psi(x) = \psi(x) - x$.
More precisely, we show that there is a direct, quantitative
connection between the size of $\Delta_\psi$ and

\begin{enumerate}

\item
the hypothetical uniformity in $\tau$,

\item
the hypothetical estimate for $F(x, T, \tau)$.

\end{enumerate}

We pick an example: see Theorem~3 and Remark~3 in
\cite{LanguascoPZ2013a} for a complete discussion of this important
topic.
Assume the RH and that $F(x, T, \tau) \ll T \log T$ uniformly for
$U \le T \le x^{1 / 2}$.
Then, if the choice
\[
  \tau
  =
  \min \Bigl( 1, \frac{\log^4 U}{\log^3 x} \Bigr),
\]
is admissible, it yields
\[
  \Delta_\psi(x)
  \ll
  x^{1 / 2} \log^2 U.
\]

\section{Asymptotic formulae}

\subsection{The asymptotic formula for
 \texorpdfstring{$F(x, T, \tau)$}{F(x, T, tau)}}

Let
\[
  S(x, \tau)
  \eqdef
  \sum_{n \ge 1} \frac{\Lambda^2(n)}n a^2(n, x, \tau),
\]
where
\[
  a(n, x, \tau)
  \eqdef
  \begin{cases}
    (n / x)^{1 / \tau} &\text{if $n \le x$} \\
    (x / n)^{1 / \tau} &\text{if $n >   x$.}
  \end{cases}
\]
Notice that $S(x, \tau) \ll \tau \log x$ by the Brun-Titchmarsh
inequality, if $\tau$ is not too small.
In \cite{LanguascoPZ2013b}, Languasco, Perelli and the author gave the
following asymptotic result.

\begin{theorem}
\label{Thm1-III}
As $x \to +\infty$ we have
\[
  F(x, T, \tau)
  \sim
  \frac T\pi \frac{S(x, \tau)}{\tau}
  +
  \frac{T \log^2 T}{\pi \tau x^{2 / \tau}}
  +
  \text{smaller order terms}
\]
uniformly for $\tau \ge 1 / T$, provided that
$T S(x, \tau) = \infty\bigl( \max (x, (\log T)^3 / \tau) \bigr)$.
\end{theorem}

Of course, this reduces to Montgomery's result for $\tau = 1$.
We remark that Theorem~\ref{Thm1-III} shows the same phenomenon
of yielding an asymptotic formula only at ``extreme ranges.''

We use a $\sigma$-uniform version of the main Lemma in Montgomery's
original paper \cite{Montgomery1973}.
He used his crucial representation for the function
\begin{equation}
\label{Mont-rep}
  (2 \sigma - 1)
  \sum_{\gamma} \frac{x^{\ii \gamma}}{(\sigma - 1 / 2)^2 + (t - \gamma)^2}
\end{equation}
for $\sigma = \frac32$, while we have to take
$\sigma = \frac12 + \tau^{-1}$.
Since we are interested in very \emph{small} values of $\tau$, our
estimates have to be uniform with respect to $\sigma$, and this leads
to a fair amount of complication in detail.
We also have to take care of other uniformity aspects.

The first ``main term'' for $F(x, T, \tau)$ arises as follows: our
version of Montgomery's representation for \eqref{Mont-rep} is the sum
of several terms, the most important being
\[
  R_1(x, t, \tau)
  \eqdef
  -x^{- 1 / 2}
  \Bigl(
    \sum_{n \le x} \Lambda(n) \Bigl( \frac xn \Bigr)^{1 - \sigma + \ii t}
    +
    \sum_{n > x}  \Lambda(n) \Bigl( \frac xn \Bigr)^{\sigma + \ii t}
  \Bigr),
\]
where, here and until the end of the section,
$\sigma = \frac12 + \tau^{-1}$.
We square out and integrate over $[- T, T]$ using Corollary~3 of
Montgomery \& Vaughan \cite{MontgomeryVaughan1974} for the
$L^2$-average of Dirichlet series over intervals.
Thus, neglecting error terms, we have
\[
  \int_{- T}^T \bigl\vert R_1(x, t, \tau) \bigr\vert^2 \, \dx t
  \sim
  2 T
  \sum_{n \ge 1} \frac{\Lambda(n)^2}n \, a(n, x, \tau)^2
  =
  2 T S(x, T).
\]
On the other hand, the integral of the square of \eqref{Mont-rep}
over the same interval is 
\[
  \sim
  4 \tau^2
  \int_{-\R}
    \Bigl\vert
      \sum_{\gamma \in [-T, T]} \frac{x^{i \gamma}}{1 + \tau^2 (t - \gamma)^2}
    \Bigr\vert^2 \, \dx t,
\]
by the Riemann-von Mangoldt formula of Theorem~\ref{R-vM}.
Squaring out, integrating term by term and computing the relevant
residue, we conclude that this is
\[
  \sim
  2 \pi \tau F(x, T, \tau).
\]
The residue turns out to be a constant multiple of
$\tau^{-1} \omega( \tau (\gamma_1 - \gamma_2) )$.

The second ``main term'' in Theorem~\ref{Thm1-III} arises computing
the contribution of another term in the development of
\eqref{Mont-rep}, namely
\[
  R_2(x, t, \tau)
  \eqdef
  \frac14
  x^{1 / 2 - \sigma + \ii t}
  \bigl( \log(\sigma^2 + t^2) + \log( (1 - \sigma)^2 + t^2) \bigr).
\]
Of course, this second main term is negligible unless $x$ is quite
small.

\subsection{The asymptotic formula for
  \texorpdfstring{$S(x, \tau)$}{S(x, tau)}}

As remarked above, if $\tau \ge x^{\eps - 1}$ then
$S(x, \tau) \ll \tau \log x$.
Moreover, if $y \le x$ and
\[
  \psi(x + y)
  -
  \psi(x)
  \sim
  y
  \qquad\text{uniformly for $y \ge x^{\beta + \eps}$}
\]
then
\[
  S(x, \tau)
  \sim
  \tau \log x
  \qquad\text{uniformly for $\tau \ge x^{\beta + \eps - 1}$.}
\]
However, $S$ is erratic for $\tau \le 1 / x$.
Essentially, it reduces to the single term given by the prime power
closest to $x$.
If $x$ itself is a prime number, then
$S(x, 1 / x) \asymp (\log x)^2 / x$, but if the prime power nearest to
$x$ is, say, $x \pm c \log x$ for some $c > 0$, then
$S(x, 1 / x) \asymp (\log x)^2 / x^{1 + c}$, which is much smaller.
Of course, gaps between consecutive primes of size $\gg \log x$ occur
infinitely often, by the PNT.

\subsection{The asymptotic formula for
  \texorpdfstring{$J(x, \tau, \theta)$}{J(x, tau, theta)}}

Finally, in \cite{LanguascoPZ2013b}, Languasco, Perelli and the author
also gave the following asymptotic result.

\begin{theorem}
\label{Thm2-III}
Assume the ``Generalised Montgomery Conjecture.''
Then
\[
  J(x, \tau, \theta)
  \sim
  \Bigl( 1 + \frac\tau2 \Bigr)
  \tau \theta x^2 \log( 1 / \theta)
\]
uniformly for $1 / x \le \theta \le x^{-\eps}$ and
$\theta^{1 / 2 - \eps} \le \tau \le 1$.
\end{theorem}

Of course, the first factor here is relevant only if $\tau \gg 1$,
when Theorem~\ref{Thm2-III} is a consequence of earlier results.
We left it for ease of comparison.
The proof requires a suitable, stronger version of the technique
introduced by Goldston \& Montgomery in \cite{GoldstonM1987}, with
particular care for the $\tau$-uniformity aspect.

\appendix
\section{The tools}

Here we briefly summarise the main tools that we need.

\begin{theorem}[Brun-Titchmarsh inequality]
\label{BT}
For $x > 0$ and $y \ge 2$ we have
\[
  \pi(x + y) - \pi(x)
  \le
  \frac{2 y}{\log y}
  \bigl( 1 + \Odi{ (\log y)^{-1} } \bigr).
\]
\end{theorem}

This is the special case $q = 1$ of Theorem~3.9 of Montgomery \&
Vaughan \cite{MontgomeryV2007}.
The error term may be omitted.
The two results that follow are classical: they are Theorem~12.5 and
Corollary 14.3 respectively of \cite{MontgomeryV2007}.

\begin{theorem}[Explicit formula]
\label{EF}
If $x \ge 2$ is not an integer and $T \ge 2$, we have
\begin{align*}
  \psi(x)
  &=
  x
  -
  \sum_{\substack{\zeta(\beta + \ii \gamma) = 0 \\ \beta \in (0, 1) \\ |\gamma| \le T}}
    \frac{x^{\rho}}{\rho}
  -\frac{\zeta'}{\zeta}(0)
  -
  \frac12
  \log \bigl(1 - x^{-2} \bigr) \\
  &\qquad\qquad+
  \Odi{\frac xT (\log xT)^2
       +
       (\log x) \min \Bigl(1, \frac x{T \langle x \rangle} \Bigr)},
\end{align*}
where $\langle x \rangle$ denotes the distance from $x$ to the nearest
prime power different from $x$ itself.
If $x$ is an integer, the term $-\frac12 \Lambda(x)$ is to be added to
the left hand side.
\end{theorem}

This is a very convenient representation of
$\Delta_\psi(x) = \psi(x) - x$, since it contains the ``free''
parameter $T$ (free in the sense that it does not appear in the
left-hand side).
This means that its value can be chosen in an optimal way in different
applications.
The upper bound in \eqref{RH} is an almost immediate deduction of the
Riemann Hypothesis, using also the Riemann-von Mangoldt formula below.

\begin{theorem}[Riemann-von Mangoldt formula]
\label{R-vM}
Let $N(T)$ denote the number of zeros of the Riemann zeta-function
lying in the rectangle in the complex plane with sides $[0, 1]$ along
the real axis and $[0, T]$ along the complex one.
Then
\[
  N(T)
  =
  \frac T{2 \pi} \Bigl( \log \frac T{2 \pi} - 1 \Bigr)
  +
  \Odi{\log T}.
\]
\end{theorem}

\begin{theorem}[Landau's formula]
\label{Landau}
We have
\[
  \Sigma(x, T)
  =
  -\frac{\Lambda(n_x)}{2 \pi}
  \frac{e^{\ii T \log(x / n_x)} - 1}{\ii \log(x / n_x)}
  +
  \Odi{x (\log(x T))^2 + \frac{\log T}{\log x}},
\]
where $n_x$ is the prime-power nearest to $x$ and $\Sigma$ is the
exponential sum over zeros defined in \textup{\eqref{def-Sigma}}.
If $x = n_x$, then the main term is $- T \Lambda(x) / 2 \pi$.
\end{theorem}

This result can be found in Ford \& Zaharescu \cite{FordZ2005} and it
is unconditional.

\begin{Conjecture}[Gonek]
\label{Gonek-C}
For $x$, $T \ge 2$ we have
\[
  \Sigma(x, T)
  \ll
  T x^{-1 / 2 + \eps} + T^{1 / 2} x^\eps.
\]
\end{Conjecture}

\vspace{0.5cm} \indent {\it
A\,c\,k\,n\,o\,w\,l\,e\,d\,g\,m\,e\,n\,t\,s.\;}
I warmly thank my coauthors Alessandro Languasco and Alberto Perelli,
and the organisers of the ``Third Italian Number Theory Meeting.''
Thanks are also due to the anonymous Referee for pointing out some
inaccuracies.

\bigskip

\phantomsection
\addcontentsline{toc}{chapter}{References}

\providecommand{\MR}{\relax\ifhmode\unskip\space\fi MR }
\providecommand{\MRhref}[2]{%
  \href{http://www.ams.org/mathscinet-getitem?mr=#1}{#2}
}
\providecommand{\href}[2]{#2}

\bigskip
\bigskip
\begin{minipage}[t]{10cm}
\begin{flushleft}
\small{
\textsc{Alessandro Zaccagnini}
\\*Dipartimento di Matematica e Informatica
\\*Universit\`a di Parma
\\*Parco Area delle Scienze, 53/a
\\*Parma, 43124 Italia
\\*e-mail: alessandro.zaccagnini@unipr.it
}
\end{flushleft}
\end{minipage}


\begin{thebibliography}{10}

\bibitem{BrudernCP1997}
\textsc{J.~Br{\"u}dern, R.~J. Cook, and A.~Perelli},
  \emph{The values of binary linear
  forms at prime arguments}, Proc. of Sieve Methods, Exponential sums and their
  Application in Number Theory (G.~R. H.~Greaves {\it et al}, ed.), Cambridge
  University Press, 1997, pp.~87--100.

\bibitem{Chan2003}
\textsc{T.~H. Chan},
  \emph{More precise pair correlation of zeros and primes in short
  intervals}, J. London Math. Soc. (2) \textbf{68} (2003), no.~3, 579--598.

\bibitem{FordZ2005}
\textsc{K.~Ford and A.~Zaharescu},
  \emph{On the distribution of imaginary parts of zeros
  of the {Riemann} zeta function}, J. reine angew. Math. \textbf{579} (2005),
  145--158.

\bibitem{Goldston2005}
\textsc{D.~A. Goldston},
  \emph{Notes on pair correlation of zeros and prime numbers},
  Recent Perspectives in Random Matrix Theory and Number Theory (F.~Mezzadri
  and N.~C. Snaith, eds.), LMS Lecture Note Series, vol. 322, Cambridge
  University Press, 2005, pp.~79--110.

\bibitem{HeathBrownG1984}
\textsc{D.~A. Goldston and D.~R. Heath-Brown},
\emph{A note on the differences between
  consecutive primes}, Math. Ann. \textbf{266} (1984), no.~3, 317--320.

\bibitem{GoldstonM1987}
\textsc{D.~A. Goldston and H.~L. Montgomery},
\emph{Pair correlation of zeros and primes
  in short intervals}, Analytic Number Theory and Diophantine Problems (Boston)
  (A.C.~Adolphson et~al., ed.), Birkh{\"a}user, 1987, pp.~183--203.

\bibitem{Gonek1993}
\textsc{S.~M. Gonek},
  \emph{An explicit formula of {Landau} and its applications to the
  theory of the zeta-function}, Contemporary Math. \textbf{143} (1993),
  395--413.

\bibitem{HeathBrown1982a}
\textsc{D.~R. Heath-Brown},
  \emph{Gaps between primes, and the pair correlation of zeros
  of the zeta-function}, Acta Arith. \textbf{41} (1982), 85--99.

\bibitem{HeathBrown1988}
\textsc{D.~R. Heath-Brown},
  \emph{The number of primes in a short interval}, J. reine angew. Math.
  \textbf{389} (1988), 22--63.

\bibitem{LanguascoPZ2012a}
\textsc{A.~Languasco, A.~Perelli, and A.~Zaccagnini},
  \emph{Explicit relations between pair correlation of zeros and primes
  in short intervals}, J. Math. Anal. Appl. \textbf{394} (2012), 761--771.

\bibitem{LanguascoPZ2013a}
\textsc{A.~Languasco, A.~Perelli, and A.~Zaccagnini},
  \emph{An extension of the pair-correlation conjecture and
  applications}, to appear in Math. Res. Lett.

\bibitem{LanguascoPZ2013b}
\textsc{A.~Languasco, A.~Perelli, and A.~Zaccagnini},
  \emph{An extended pair-correlation conjecture and primes in short
  intervals}, to appear in Trans. Amer. Math. Soc.

\bibitem{Littlewood1914}
\textsc{J.~E. Littlewood},
  \emph{Sur la distribution des nombres premiers}, C. R. Acad.
  Sc. Paris \textbf{158} (1914), 1869--1872.

\bibitem{Maier1985}
\textsc{H.~Maier},
  \emph{Primes in short intervals}, Michigan Math. J. \textbf{32}
  (1985), 221--225.

\bibitem{Montgomery1973}
\textsc{H.~L. Montgomery},
  \emph{The pair correlation of zeros of the zeta function},
  Analytic number theory ({Proc}. {Sympos}. {Pure} {Math}., {Vol}. {XXIV},
  {St}. {Louis} {Univ}., {St}. {Louis}, {Mo}., 1972), Amer. Math. Soc.,
  Providence, R.I., 1973, pp.~181--193.

\bibitem{MontgomeryS2002}
\textsc{H.~L. Montgomery and K.~Soundararajan},
  \emph{Beyond pair correlation}, Paul
  {Erd{\H o}s} and his mathematics, I (Budapest, 1999), Bolyai Soc. Math.
  Stud., vol.~11, J\'anos Bolyai Math. Soc., 2002, pp.~507--514.

\bibitem{MontgomeryVaughan1974}
\textsc{H.~L. Montgomery and R.~C. Vaughan},
  \emph{Hilbert's inequality}, J. London
  Math. Soc. (2) \textbf{8} (1974), 73--82.

\bibitem{MontgomeryV2007}
\textsc{H.~L. Montgomery and R.~C. Vaughan},
  \emph{Multiplicative {Number} {Theory}. {I}. {Classical} {Theory}},
  Cambridge University Press, Cambridge, 2007.

\bibitem{MurtyP1999}
\textsc{M.~R. Murty and A. Perelli},
  \emph{The pair correlation of zeros of functions in the {Selberg} class},
  Internat. Math. Res. Notices (1999), no.~10, 531--545.

\bibitem{Pintz1984}
\textsc{J.~Pintz},
  \emph{On the remainder term of the prime number formula and the zeros
  of the {Riemann} zeta-function}, Number Theory (Noordwijkerhout), Lecture
  Notes in Mathematics, vol. 1068, Springer, 1984, pp.~186--197.

\bibitem{Riemann1859}
\textsc{G.~F.~B. Riemann},
  \emph{{\"Uber} die {Anzahl} der {Primzahlen} unter einer
  gegebenen {Gr\"osse}}, Monatsber. {K\"onigl.} {Preuss.} {Akad.} {Wiss.}
  Berlin (1859), 671--680, in ``Gesammelte Mathematische Werke'' (ed. H.
  Weber), Dover reprint 1953.

\bibitem{SaffariV1977}
\textsc{B.~Saffari and R.~C. Vaughan},
  \emph{On the fractional parts of $x/n$ and
  related sequences. {II}}, Ann. Inst. Fourier \textbf{27} (1977), 1--30.

\bibitem{Selberg1943}
\textsc{A.~Selberg},
  \emph{On the normal density of primes in small intervals, and the
  difference between consecutive primes}, Arch. Math. Naturvid. \textbf{47}
  (1943), 87--105.

\bibitem{Zaccagnini1998}
\textsc{A.~Zaccagnini},
  \emph{Primes in almost all short intervals}, Acta Arith.
  \textbf{84.3} (1998), 225--244.

\bibitem{Zaccagnini2000b}
\textsc{A.~Zaccagnini},
  \emph{A conditional density theorem for the zeros of the {Riemann}
  zeta-function}, Acta Arith. \textbf{93} (2000), 293--301.

\end{thebibliography}
\end{document}